\documentclass[12pt]{article}
\usepackage{amsfonts}
\usepackage{amssymb}
\usepackage{amsthm}
\usepackage[centertags]{amsmath}
\usepackage{newlfont}
\usepackage{graphicx}

\newfont{\bb}{msbm10 at 12pt}
\def\r{\hbox{\bb R}}

\def\e{\hbox{\bf E}}
\def\t{\hbox{\bf T}}
\def\n{\hbox{\bf N}}
\def\b{\hbox{\bf B}}

\newtheorem{theorem}{Theorem}[section]
\newtheorem{definition}[theorem]{Definition}

\newtheorem{remark}[theorem]{Remark}

\begin{document}

\title{On slant helices in Minkowski space $\e_1^3$}
\author{ Ahmad T. Ali\\Mathematics Department\\
 Faculty of Science, Al-Azhar University\\
 Nasr City, 11448, Cairo, Egypt\\
email: atali71@yahoo.com\\
\vspace*{1cm}\\
 Rafael L\'opez\footnote{Partially
supported by MEC-FEDER
 grant no. MTM2007-61775.}\\
Departamento de Geometr\'{\i}a y Topolog\'{\i}a\\
Universidad de Granada\\
18071 Granada, Spain\\
email: rcamino@ugr.es}
\date{}

\maketitle
\begin{abstract} We consider a curve $\alpha=\alpha(s)$ in Minkowski 3-space $\e_1^3$ and denote by $\{\t,\n,\b\}$ the Frenet frame of $\alpha$. We say that $\alpha$ is a slant helix if there exists a fixed direction $U$ of $\e_1^3$ such that the function $\langle \n(s),U\rangle$ is constant.   In this work we give characterizations of slant helices in terms of the curvature and torsion of $\alpha$.
\end{abstract}



%%%%%%%%%%%%%%%%%%%%%%%%%%%%%%%%%%%%%%%%%%%%%%%
\section{Introduction and statement of results}
%%%%%%%%%%%%%%%%%%%%%%%%%%%%%%%%%%%%%%%%%%%%%%%

Let $\e_1^3$ be the   Minkowski 3-space, that is,  $\e_1^3$ is the real vector space $\r^3$ endowed with the standard flat metric
$$\langle,\rangle=dx_1^2+dx_2^2-dx_3^2,$$
where $(x_1,x_2,x_3)$ is a rectangular coordinate system of $\e_1^3$. An arbitrary vector $v\in\e_1^3$ is said spacelike if $\langle v,v\rangle>0$ or $v=0$, timelike if $\langle v,v\rangle<0$, and lightlike (or null) if $\langle v,v\rangle =0$ and $v\neq0$. The norm (length) of a vector $v$ is given by $\parallel v\parallel=\sqrt{|\langle v,v\rangle|}$.

Given a regular (smooth) curve $\alpha:I\subset\r\rightarrow\e_1^3$, we say that $\alpha$ is spacelike (resp.  timelike,  lightlike) if all of its velocity vectors $\alpha'(t)$ are  spacelike (resp. timelike, lightlike). If $\alpha$ is spacelike or timelike we say that $\alpha$ is a non-null curve. In such case, there exists a change of the parameter $t$, namely, $s=s(t)$, such that $\parallel\alpha'(s)\parallel=1$.  We say then that $\alpha$ is parametrized by the arc-length parameter.  If the curve $\alpha$ is lightlike, the acceleration vector $\alpha''(t)$ must be spacelike for all $t$. Then we change the parameter $t$ by $s=s(t)$ in such way that $\parallel \alpha''(s)\parallel=1$ and we say that $\alpha$ is parameterized by the pseudo arc-length parameter. In any of the above cases, we say that $\alpha$ is a unit speed curve.

Given a unit speed curve $\alpha$ in Minkowski space $\e_1^3$ it is possible to define a Frenet frame $\{\t(s),\n(s),\b(s)\}$ associated for each point $s$ \cite{ku,wa}. Here $\t$, $\n$ and $\b$ are the tangent, normal and binormal vector field, respectively. The geometry of the curve $\alpha$ can be describe by the differentiation of the Frenet frame, which leads to the corresponding Frenet equations. Although different expressions of the Frenet equations appear depending of the causal character of the Frenet trihedron (see the next sections below), we have the concepts of curvature $\kappa$ and torsion $\tau$ of the curve. With this preparatory introduction, we give the following

\begin{definition} A unit speed curve $\alpha$ is called a slant helix if there exists a constant vector field $U$ in $\e_1^3$ such that the function $\langle \n(s),U\rangle$ is constant.
\end{definition}

This definition is motivated by what happens in Euclidean ambient space $\e^3$. In this setting, we recall that a helix is a curve where the tangent lines make a constant angle with a fixed direction. Helices are characterized by the fact that the ratio $\tau/\kappa$ is constant along the curve \cite{dc}. Helices in Minkowski space have been studied depending on the causal character of the curve $\alpha$: see for example \cite{fgl,ko,ps}.
 Recently, Izumiya and Takeuchi have introduced the concept of slant helix in Euclidean space by saying that the normal lines make a constant angle with a fixed direction \cite{it}. They characterize a slant helix if and only if the function
\begin{equation}\label{slant}
\dfrac{\kappa^2}{(\kappa^2+\tau^2)^{3/2}}\Big(\dfrac{\tau}{\kappa}\Big)'
\end{equation}
is constant. See also \cite{ky, okkk}. Thus, our definition of slant helix is the Lorentzian version of the Euclidean one. Only it is important to point out that, in contrast to what happens in Euclidean space, in Minkowski ambient space we can not  define the angle between two vectors (except that both vectors are of  timelike type). For this reason, we avoid to say about the angle between the vector fields $\n(s)$ and $U$.

Our main result in this work is the following characterization of slant helices in the spirit of the one given in equation (\ref{slant}). We will assume throughout this work that the curvature and torsion functions do not  equal zero. Exactly, we prove

\begin{theorem}\label{t1} Let $\alpha$ be a unit speed timelike curve in $\e_1^3$. Then $\alpha$ is a slant helix if and only if either one the next two functions
\begin{equation}\label{slant2}
\frac{\kappa^2}{(\tau^2-\kappa^2)^{3/2}}\Big(\dfrac{\tau}{\kappa}\Big)'\hspace*{1cm}\mbox{or}\hspace*{1cm}
\frac{\kappa^2}{(\kappa^2-\tau^2)^{3/2}}\Big(\dfrac{\tau}{\kappa}\Big)'
\end{equation}
is constant everywhere $\tau^2-\kappa^2$ does not vanish.
\end{theorem}

\begin{theorem}\label{t2} Let $\alpha$ be a unit speed spacelike curve in $\e_1^3$.
\begin{enumerate}
\item If the normal vector of $\alpha$ is spacelike, then $\alpha$ is a slant helix if and only if either one the next two functions
\begin{equation}\label{slant3}
\frac{\kappa^2}{(\tau^2-\kappa^2)^{3/2}}\Big(\dfrac{\tau}{\kappa}\Big)'\hspace*{1cm}\mbox{or}\hspace*{1cm}
\frac{\kappa^2}{(\kappa^2-\tau^2)^{3/2}}\Big(\dfrac{\tau}{\kappa}\Big)'
\end{equation}
is constant everywhere $\tau^2-\kappa^2$ does not vanish.
\item If the normal vector of $\alpha$ is timelike, then $\alpha$ is a slant helix if and only if the function
\begin{equation}\label{slant4}
\frac{\kappa^2}{(\tau^2+\kappa^2)^{3/2}}\Big(\dfrac{\tau}{\kappa}\Big)'
\end{equation}
is constant.
\item Any spacelike curve with lightlike normal vector is a slant curve. 
\end{enumerate}
\end{theorem}

In the case that $\alpha$ is a lightlike curve, we have

\begin{theorem}\label{t3} Let $\alpha$ be a unit speed lightlike curve in $\e_1^3$. Then $\alpha$ is a slant helix if and only if the torsion is
\begin{equation}\label{slant5}
\tau(s)=\frac{a}{(bs+c)^2},
\end{equation}
where $a$, $b$ and $c$ are constant.
\end{theorem}

The proof of Theorems \ref{t1}, \ref{t2} and \ref{t3}  is carried in the successive sections.

%%%%%%%%%%%%%%%%%%%%%%%%%%%%%%%%%%%%%%%%%%%%%%%%%%%%%%%%%%
\section{Timelike  slant helices}
%%%%%%%%%%%%%%%%%%%%%%%%%%%%%%%%%%%%%%%%%%%%%%%%%%%%%%%%%%%%
Let $\alpha$ be a unit speed timelike curve in $\e_1^3$. The Frenet frame $\{\t,\n,\b\}$ of $\alpha$ is given by
$$\t(s)=\alpha'(s),\ \ \n(s)=\dfrac{\alpha''(s)}{\parallel\alpha''(s)\parallel},\ \ \b(s)=\t(s)\times\n(s).$$
The  Frenet equations are
\begin{equation}\label{equi1}
 \left[
   \begin{array}{c}
     \t'(s) \\
     \n'(s) \\
     \b'(s)
        \end{array}
 \right]=\left[
           \begin{array}{ccc}
             0 & \kappa(s) & 0 \\
             \kappa(s) & 0 &\tau(s)\\
             0 &-\tau(s) & 0\\
           \end{array}
         \right]\left[
   \begin{array}{c}
     \t(s) \\
     \n(s) \\
     \b(s) \\
   \end{array}
 \right].
 \end{equation}
 In order to prove Theorem \ref{t1}, we first assume that $\alpha$ is a slant helix. Let $U$ be the vector field  such that the function $\langle \n(s),U\rangle:=c$ is constant. There exist smooth functions $a_1$ and $a_3$ such that
\begin{equation}\label{u1}
U=a_1(s)\t(s)+c \n(s)+a_3(s) \b(s),\ \ s\in I.
\end{equation}
As  $U$ is constant, a differentiation in (\ref{u1}) together (\ref{equi1}) gives
\begin{equation}\label{u2}
\left.\begin{array}{ll}
a_1'-c\kappa&=0\\
\kappa a_1-\tau a_3&=0\\
a_3'+c \tau &=0
\end{array}\right\}
\end{equation}
From the second equation in (\ref{u2}) we have
\begin{equation}\label{u5}
a_1=a_3\big(\dfrac{\tau}{\kappa}\big).
\end{equation}
Moreover
\begin{equation}\label{u3}
\langle U,U\rangle=-a_1^2+c^2+a_3^2=\mbox{constant}.
\end{equation}
 We point out that this constraint, together the second and third equation of (\ref{u2}) is equivalent to the very system (\ref{u2}). From (\ref{u5}) and (\ref{u3}), set
$$a_3^2\Big(\big(\frac{\tau}{\kappa}\big)^2-1\Big)=\epsilon m^2,\ \ m>0,\epsilon\in\{-1,0,1\}.$$
If $\epsilon=0$, then $a_3=0$ and from (\ref{u2}) we have $a_1=c=0$. This means that $U=0$: contradiction. Thus $\epsilon=1$ or $\epsilon=-1$ which gives
$$a_3=\pm\dfrac{m}{\sqrt{\big(\dfrac{\tau}{\kappa}\big)^2-1}}\hspace*{1cm}\mbox{or}\hspace*{1cm}a_3=\pm\dfrac{m}{\sqrt{1-\big(\dfrac{\tau}{\kappa}\big)^2}}$$
on $I$. The third equation in (\ref{u2}) yields
$$\dfrac{d}{ds}\Big[\pm\dfrac{m}{\sqrt{\big(\dfrac{\tau}{\kappa}\big)^2-1}}\Big]=-c \tau \hspace*{1cm} \mbox{or}\hspace*{1cm}\dfrac{d}{ds}\Big[\pm\dfrac{m}{\sqrt{1-\big(\dfrac{\tau}{\kappa}\big)^2}}\Big]=c \tau$$
on $I$. This can be written as
$$\frac{\kappa^2}{(\tau^2-\kappa^2)^{3/2}}\Big(\dfrac{\tau}{\kappa}\Big)'=\mp\dfrac{c}{m}\hspace*{1cm}\mbox{or}\hspace*{1cm}
\frac{\kappa^2}{(\kappa^2-\tau^2)^{3/2}}\Big(\dfrac{\tau}{\kappa}\Big)'=\pm\dfrac{c}{m}$$
This shows a part of Theorem \ref{t1}. Conversely, assume that the condition (\ref{slant2}) is satisfied. In order to simplify the computations, we assume that the first function in (\ref{slant2}) is a constant, namely,  $c$ (the other case is analogous). We define
\begin{equation}\label{u9}
U=\dfrac{\tau}{\sqrt{\tau^2-\kappa^2}}\t+
c\n+\dfrac{\kappa}{\sqrt{\tau^2-\kappa^2}}\b\Big.
\end{equation}
A differentiation of (\ref{u9}) together the Frenet equations gives $\dfrac{dU}{ds}=0$, that is,  $U$ is a constant vector.
On the other hand, $\langle\n(s),U\rangle=1$ and this means that $\alpha$ is a slant helix.

\begin{remark}  In Theorem \ref{t1} we need to assure that the function $\tau^2-\kappa^2$ does not vanish everywhere. We do not know that happens if it vanishes at some points. On the other hand, any timelike curve that satisfies $\tau(s)^2-\kappa(s)^2=0$ is a slant curve. The reasoning is the following. For simplicity, we only consider  the case  that $\tau=\kappa$. We define $U=\t(s)+\b(s)$, which is constant using the Frenet equations (\ref{equi1}). Moreover,
 $\langle \n,U\rangle=0$, that is, $\alpha$ is a slant curve. Finally, we point that there exist curves in $\e_1^3$ that satisfies the relation
 $\tau=\kappa$: it suffices to put $\tau=\kappa:=c=\mbox{constant}$ and the fundamental theorem of the theory of curves assures the existence of a timelike curve $\alpha$ with curvature and torsion $c$.
\end{remark}
%%%%%%%%%%%%%%%%%%%%%%%%%%%%%%%%%%%%%%%%%%%%%%%%%%%%%%%%%%
\section{Spacelike  slant helices}
%%%%%%%%%%%%%%%%%%%%%%%%%%%%%%%%%%%%%%%%%%%%%%%%%%%%%%%%%%%%

Let $\alpha$ be a unit speed spacelike curve in $\e_1^3$. In the case that
the normal vector $\n(s)$ of $\alpha$ is spacelike or timelike, the proof of Theorem \ref{t2} is similar to the one given for Theorem \ref{t1}. We omit the details.

The  case that remains to study is that the normal vector $\n(s)$ of the curve is a lightlike vector for any $s\in I$. Now the Frenet trihedron is
$\t(s)=\alpha'(s)$, $\n(s)=\t'(s)$ and $\b(s)$ is the unique lightlike vector orthogonal to $\t(s)$ such that $\langle\n(s),\b(s)\rangle=1$. Then the Frenet equations as
\begin{equation}\label{u11}
 \left[
   \begin{array}{c}
     \t' \\
     \n' \\
     \b'
        \end{array}
 \right]=\left[
           \begin{array}{ccc}
             0 & 1 & 0 \\
             0 & \tau & 0 \\
             -1 & 0 & \tau\\
           \end{array}
         \right]\left[
   \begin{array}{c}
     \t \\
     \n \\
     \b \\
   \end{array}
 \right].
 \end{equation}
Here $\tau$ is the torsion of the curve (recall that $\tau(s)\not=0$ for any $s\in I$). We show that \emph{any} such curve is a slant helix. Let $a_2(s)$ any non-trivial solution of the O.D.E. $y'(s)+\tau(s)y(s)=0$ and define $U=a_2(s)\n(s)$. By using (\ref{u11}),
$dU(s)/ds=0$, that is, $U$ is a (non-zero) constant vector field of $\e_1^3$ and, obviously, the function $\langle\n(s),U\rangle$ in constant (and equal to $0$).

%%%%%%%%%%%%%%%%%%%%%%%%%%%%%%%%%%%%%%%%%%%%%%%%%%%%%%%%%%
\section{Lightlike slant helices }
%%%%%%%%%%%%%%%%%%%%%%%%%%%%%%%%%%%%%%%%%%%%%%%%%%%%%%%%%%%%

In this section we show Theorem \ref{t3}. Let $\alpha$ be a unit lightlike in $\e_1^3$.  The Frenet frame of $\alpha$ is
$\t(s)=\alpha'(s)$, $\n(s)=\t'(s)$ and $\b(s)$ the unique lightlike vector orthogonal to $\n(s)$ such that
$\langle\t(s),\b(s)\rangle=1$. The Frenet equations are
\begin{equation}\label{u21}
 \left[
   \begin{array}{c}
     \t' \\
     \n' \\
     \b'
        \end{array}
 \right]=\left[
           \begin{array}{ccc}
             0 & 1 & 0 \\
             \tau & 0 & -1 \\
             0 & -\tau & 0\\
           \end{array}
         \right]\left[
   \begin{array}{c}
     \t \\
     \n \\
     \b \\
   \end{array}
 \right].
 \end{equation}
Here $\tau(s)$ is the torsion of $\alpha$, which is assumed with the property $\tau(s)\not=0$, for any $s\in I$.

Assume that $\alpha$ is a slant helix. Let $U$  be the constant vector field such that the function $\langle \n(s),U\rangle$ is constant.
As in the above cases
$$U=a_1(s)\t(s)+c \n(s)+a_3(s) \b(s),\ \ s\in I,$$
where $c$ is a constant and
\begin{equation}\label{u23}
\left.\begin{array}{ll}
a_1'+c \tau&=0\\
a_1-\tau a_3&=0\\
a_3'-c &=0
\end{array}\right\}
\end{equation}
Then $a_3(s)=cs+m$, $m\in \r$ and $a_1=(cs+m)\tau$. Using the first equation of (\ref{u23}), we have
$(cs+m)\tau'+2c\tau=0$. The solution of this equation is
$$\tau(s)=\frac{n}{(cs+m)^2},$$
where $m$ and $n$ are constant. This proves (\ref{slant5}) in Theorem \ref{t3}. Conversely, if the condition (\ref{slant5}) is satisfied,
 we define
 $$U=\frac{a}{bs+c}\t(s)+b\n(s)+(bs+c)\b(s).$$
Using the Frenet equations (\ref{u21}) we obtain that $dU(s)/ds=0$, that is, $U$ is a constant vector field of $\e_1^3$. Finally, $\langle \n(s),U\rangle=b$ and this proves that $\alpha$ is a slant helix.

\end{document}